\documentclass{statsoc}

\def\eqref#1{\mbox{(\ref{eq:#1})}}    

\title{On the Generalized Poisson Distribution}
\author{Hans J. H. Tuenter}
\address{Schulich School of Business, 
         York University, Toronto, 
         Canada, M3J 1P3}
\email{HTuenter@Schulich.YorkU.Ca}

\begin{document}

\begin{abstract}
We use Euler's difference lemma
to prove that,
for $\theta>0$ and $0\le\lambda<1$, 
the function $P_n$ defined on the non-negative integers by
\[ P_n(\theta,\lambda)={\theta(\theta+n\lambda)^{n-1}\over n!}
         e^{-n\lambda-\theta}
\]
defines a probability distribution,
known as the Generalized Poisson Distribution.
\end{abstract}

\section{Introduction}
The Generalized Poisson Distribution (GPD), introduced in~\cite{ConsulJain73},
and studied extensively by~\cite{Consul89}
is defined on the non-negative integers,
for $0\le\lambda<1$ and $\theta>0$, by
\begin{equation}
  P_n(\theta,\lambda)={\theta(\theta+n\lambda)^{n-1}\over n!}
       e^{-\theta-n\lambda}.
  \label{eq:1}
\end{equation}

Applications of the GPD can be found in settings 
where one seeks to describe the distribution of an event that occurs 
rarely in a short period, but where we observe the frequency of its occurrence 
in longer periods of time.
It extends the Poisson distribution by its ability to describe situations
where the probability of occurrence of a single event does not remain constant (as in a Poisson process), but is affected by previous occurrences.
The distribution has been found~\cite[pp. 117--129]{Consul89} to accurately describe phenomena as diverse as the observed number of industrial accidents and injuries, where a learning effect may be present,
the spatial distribution of insects, where initial occupation of a spot by a member of the species has an influence on the attractiveness of the spot to other members of the species,
and the number of units of different commodities purchased by consumers, where current sales have an impact on the level of subsequent sales through repeat purchases.

It was shown by~\cite{ConsulJain73} that~\eqref{1} is a probability distribution, as it has the property $\sum_{n=0}^{\infty} P_n(\theta,\lambda)=1$, 
by using an identity that can be found in~\citet[eq.~6]{Jensen02}.
However, as remarked by~\citet[p.~12]{Consul89},
``It is very difficult to prove by direct summation that the sum
of all the probabilities is unity''.
Recently,~\cite{LernerLR97} gave a more direct proof
using analytic functions.
We give a shorter and more elegant proof
based upon an application of Euler's classic difference lemma.

\section{Derivation}
We shall first prove the identity
\begin{equation}
 \sum_{n=0}^{\infty} {(\theta+\lambda n)^n \over n!} e^{-\theta-\lambda n} = 
 {1\over1-\lambda}, 
 \mbox{\rm \ \ for\ } 
 -\lambda_0<\lambda<1,
 \label{eq:2}
\end{equation}
where $\lambda_0=0.2784645428\ldots$ is the solution to $\lambda e^{\lambda}=e^{-1}$.
Let $S(\theta,\lambda)$ denote the sum in~\eqref{2}.
Expanding the exponential, grouping terms and changing the order of summation gives 
\begin{eqnarray}
  S(\theta,\lambda)
  & = &
  \sum_{n=0}^{\infty} {(\theta+\lambda n)^n\over n!} 
  \sum_{k=0}^{\infty} (-1)^k  {(\theta+\lambda n)^k\over k!}  
   \ = \ 
  \sum_{n=0}^{\infty} \sum_{k=n}^{\infty} 
     (-1)^{k-n}
   {(\theta+\lambda n)^k\over n!\ (k-n)!}  \label{eq:3} \\
  & = &
 \sum_{k=0}^{\infty} {1\over k!}  
 \sum_{n=0}^k  (-1)^{k-n} \left({k \atop n}\right)
   (\theta+\lambda n)^k . 
 \label{eq:4}
\end{eqnarray}
Now use Euler's difference formula,
which states that the $k$th difference of any $p$th power is
$0$ for $p<k$ and $k!$ times the leading coefficient of the $p$th
power for $p=k$. 
The most convenient representation for our purpose, as can be found 
in~\citet[eq.~5.12]{Gould78}, is 
\begin{eqnarray*}
  \sum_{n=0}^k
  (-1)^{k-n}
  \left({k \atop n}\right)
  (A+Bn)^p
  & = &
  \left\{%
  \begin{tabular}{l l}
  $0$,      & $0\le p <k$  \\
  $B^k k!$, & $p=k$
  \end{tabular}
  \right..  
  \label{eq:5}
\end{eqnarray*}
Applying this to the inner summation in~\eqref{4} gives 
$ S(\theta,\lambda)
  =
  \sum_{k=0}^{\infty} {1\over k!} \lambda^k k! 
  = \sum_{k=0}^{\infty} \lambda^k ={1\over 1-\lambda},
$
and establishes the identity in~\eqref{2}.

We now show that for $\left|\lambda\right|<\lambda_0$
the interchanging of the summation signs in the proof of~\eqref{2} is allowed,
as we have tacitly assumed,
by establishing absolute convergence.
Taking absolute values of the summands one sees that the inner summation in~\eqref{3} reduces to $e^{\left|\theta+\lambda n\right|}$.
Now apply Cauchy's root test and Stirling's approximation $n!\sim n^n e^{-n} \sqrt{2\pi n}$ to give
\[ \limsup_{n\rightarrow\infty} \sqrt[n]{\frac{\left|\theta+\lambda    n\right|^n}{n!}e^{\left|\theta+\lambda n\right|}}
   =\left|\lambda\right| e^{1+\left|\lambda\right|}<1,
\]
or the desired $\left|\lambda\right|<\lambda_0$ as a criterion for absolute convergence.
Another application of the root test shows that the left-hand side of~\eqref{2} 
converges for $\left|\lambda \smash{e^{-\lambda}}\right|<e^{-1}$,
and as all the summands are positive, this convergence is uniform.
Thus we can extend the range by analyticity to $-\lambda_0<\lambda<1$,
and this completes the proof. 
 
Using formula~\eqref{2}, we have
\begin{eqnarray*}
  \sum_{n=0}^{\infty} P_n(\theta,\lambda) 
  &=& \sum_{n=0}^{\infty}\frac{(\theta+n\lambda)^n}{n!}e^{-\theta-n\lambda} -
      \lambda\sum_{n=1}^{\infty}\frac{(\theta+n\lambda)^{(n-1)}}{(n-1)!}e^{-\theta-n\lambda} \\
  &=& S(\theta,\lambda) -
      \lambda\sum_{n=0}^{\infty}\frac{(\theta+\lambda+n\lambda)^{n}}{n!}e^{-\theta-\lambda-n\lambda} \\
  &=& S(\theta,\lambda) - \lambda S(\theta+\lambda,\lambda) = 1, 
\end{eqnarray*}
and proves that~\eqref{1} indeed defines a probability distribution.

\section{Acknowledgments}
I would like to thank one of the referees for suggesting a more concise version of the last part of the proof that $\sum_n P_n(\theta,\lambda)=1$.
Furthermore, as was pointed out by the other referee, 
it should be noted that formula~\eqref{2} can also be derived from an identity,
posed as ``problem~214'' by~\citet[pp.~126, 302]{PolyaSzego:1925}.
However, the proof in the present paper is more direct, and does not rely
on an application of Lagrange's formula. 


\end{document}